\documentclass[a4paper,10pt,reqno]{amsart}
\usepackage[a4paper,tmargin=40mm,bmargin=35mm,hmargin=30mm]{geometry}
\usepackage[T1]{fontenc}
\usepackage{lmodern}
\usepackage{amsmath}
\usepackage{amssymb}
\usepackage{amsthm}
\usepackage{stmaryrd}
\usepackage{tikz}
\usepackage[all]{xy}
\usepackage{booktabs}
\usepackage{multirow}
\usepackage{multicol}
\usepackage{here}
\usepackage{aliascnt}
\usepackage{hyperref}
\usepackage[noabbrev]{cleveref}

\newcommand{\NewTheorem}[2]{
	\newaliascnt{#1}{TheoremEnvironment}
	\newtheorem{#1}[#1]{#1}
	\aliascntresetthe{#1}
	\crefname{#1}{#1}{#2}
	\Crefname{#1}{#1}{#2}
}

\theoremstyle{definition}

\NewTheorem{Definition}{Definitions}
\NewTheorem{Remark}{Remarks}
\NewTheorem{Axiom}{Axioms}
\NewTheorem{Example}{Examples}
\NewTheorem{Observation}{Observations}
\NewTheorem{Convention}{Conventions}
\NewTheorem{Notation}{Notations}
\NewTheorem{Setting}{Settings}
\NewTheorem{Question}{Questions}
\NewTheorem{Answer}{Answers}
\NewTheorem{Conjecture}{Conjectures}
\NewTheorem{Problem}{Problems}
\NewTheorem{Solution}{Solutions}
\NewTheorem{Goal}{Goals}
\NewTheorem{Comment}{Comments}
\NewTheorem{Aim}{Aims}
\NewTheorem{Caution}{Cautions}
\NewTheorem{Exercise}{Exercises}
\NewTheorem{Examples}{Examples}
\theoremstyle{plain}
\NewTheorem{Proposition}{Propositions}
\NewTheorem{Lemma}{Lemmas}
\NewTheorem{Theorem}{Theorems}
\NewTheorem{Corollary}{Corollaries}

\crefname{enumi}{}{}
\Crefname{enumi}{}{}
\creflabelformat{enumi}{(#2#1#3)}
\crefname{enumii}{}{}
\Crefname{enumii}{}{}
\creflabelformat{enumii}{(#2#1#3)}
\crefname{enumiii}{}{}
\Crefname{enumiii}{}{}
\creflabelformat{enumiii}{(#2#1#3)}

\makeatletter
\renewcommand{\p@enumii}{}
\renewcommand{\p@enumiii}{}
\makeatother

\numberwithin{equation}{section}
\crefname{equation}{}{}
\Crefname{equation}{}{}
\creflabelformat{equation}{(#2#1#3)}

\newcommand{\SwapSymbols}[1]{
	\expandafter\let\expandafter\temporarysymbol\csname #1\endcsname
	\expandafter\let\csname #1\expandafter\endcsname\csname var#1\endcsname
	\expandafter\let\csname var#1\endcsname\temporarysymbol
}

\SwapSymbols{epsilon}
\SwapSymbols{phi}
\SwapSymbols{Gamma}
\SwapSymbols{Delta}
\SwapSymbols{Theta}
\SwapSymbols{Lambda}
\SwapSymbols{Xi}
\SwapSymbols{Pi}
\SwapSymbols{Sigma}
\SwapSymbols{Upsilon}
\SwapSymbols{Phi}
\SwapSymbols{Psi}
\SwapSymbols{Omega}

\newcommand{\cA}{\mathcal{A}}

\newcommand{\cC}{\mathcal{C}}
\newcommand{\cD}{\mathcal{D}}

\newcommand{\cF}{\mathcal{F}}
\newcommand{\cG}{\mathcal{G}}

\newcommand{\cN}{\mathcal{N}}

\newcommand{\cS}{\mathcal{S}}
\newcommand{\cT}{\mathcal{T}}
\newcommand{\cU}{\mathcal{U}}

\newcommand{\cX}{\mathcal{X}}

\newcommand{\To}{\longrightarrow}

\DeclareMathOperator{\Hom}{Hom}

\DeclareMathOperator{\Ext}{Ext}
\DeclareMathOperator{\Tor}{Tor}

\DeclareMathOperator{\Ann}{Ann}

\DeclareMathOperator{\Ker}{Ker}
\DeclareMathOperator{\Coker}{Coker}
\let\Im\relax
\DeclareMathOperator{\Im}{Im}

\DeclareMathOperator{\Spec}{Spec}

\DeclareMathOperator{\Min}{Min}

\DeclareMathOperator{\Supp}{Supp}

\title{Cofiniteness with respect to extension of Serre subcategories at small dimensions}

\author{Reza Sazeedeh}
\address{Department of Mathematics, Urmia University, P.O.Box: 165, Urmia, Iran}
\email{rsazeedeh@ipm.ir}

\begin{document}

\begin{abstract}
Let $R$ be a commutative noetherian ring, $\frak a$ be an ideal of $R$, $\cS$ be an arbitrary Serre subcategory of $R$-modules and let $\cN$ be the subcategory of finitely generated $R$-modules. In this paper, we study $\cN\cS$-$\frak a$-cofinite modules with respect to the extension subcategory $\cN\cS$ when $\dim R/\frak a\leq 2$. We also study $\frak a$-cofiniteness with respect to a new dimension. 

\end{abstract}

\maketitle
\tableofcontents

\section{Introduction}
Throughout this paper $R$ is a commutative noetherian ring, $\frak a$ is an ideal of $R$. Given a Serre subcategory $\cS$  of $R$-modules, an $R$-module $M$ is said to be $\cS$-$\frak a$-{\it cofinite} if $\Supp M\subseteq V(\frak a)$ and $\Ext_R^i(R/\frak a, M)\in\cS$ for all integers $i\geq 0$. Let $\cN$ be the subcategory of finitely generated $R$-modules.

 The extension subcategory induced by $\cN$ and $\cS$ is denoted by $\cN\cS$, consisting of those $R$-modules $M$ for which there exist an exact sequence $0\To N\To M\To S\To 0$ such that $N\in\cN$ and $S\in\cS$. It has been proved by [Y] that $\cN\cS$ is Serre. A well-known example for these subcategories is $\cN\cA$, the subcategory of {\it minimax modules} studied by [Z] where $\cA$ is the subcategory of artinian modules. Another example is $\cN\cF$, the subcategory of FSF modules introduced by [Q], where $F$ consists of all modules of finite support. When $\cS=0$, an $\cN\cS$-$\frak a$-cofinite module was known as an $\frak a$-cofinite module which is defined for the first time by Hartshorne [H], giving a negative answer to a question of [G, Expos XIII, Conjecture 1.1]. Many people [BNS, M1, M2, NS] studied $\frak a$-cofiniteness in various cases. 

The main aim of this paper is to extend the fundamental results about $\frak a$-cofinite modules at small dimensions to $\cN\cS$-$\frak a$ cofinite modules.
 We recall that a Serre subcategory $\cS$ satisfies the condition $C_{\frak a}$ if for every $R$-module $M$, the following implication holds.
\begin{center}
$C_{\frak a}$: If $\Gamma_{\frak a}(M)=M$ and $(0:_M{\frak a})$ is in $\cS$, then
$M$ is in $\cS$.
\end{center}

For every $\frak p\in\Spec R$, we denote by $\cS({\frak p})$ the smallest Serre subcategory of $R_{\frak p}$-modules containing $\cS_{\frak p}=\{M_{\frak p}|\hspace{0.1cm} M\in\cS\}$. For every $R$-module $M$, $\dim M$ means the dimension of $\Supp_RM$ which is the length of the longest chain of prime ideals in $\Supp_RM$. In Section 2, we prove the following result. 
\begin{Theorem}
Let $(R,\frak m)$ be a local ring, let $S(\frak p)$ satisfy the condition $\frak pR_{\frak p}$ for each $\frak p\in V(\frak a)$. If $M$ is an $R$-module of dimension $d$ such that $\Ext_R^i(R/\frak a,M)\in\cN\cS$ for each $i\leq d$, then $\Ext_R^i(N,M)\in\cN\cS$ for each $i\geq 0$ and each finitely generated $R$-module $N$ with $\dim N\leq 2$ and $\Supp_RN\subseteq V(\frak a)$.
\end{Theorem}
 Assume that $M$ is an $R$-module such that $\Supp_R M\subseteq V(\frak a)$. Melkersson [M1, Theorem 2.3] showed that if $\dim R/\frak a=1$, then $M$ is $\frak a$-cofinite if and only if $\Hom_R(R/\frak a,M)$ and $\Ext_R^1(R/\frak a,M)$ are finitely generated. The above theorem generalizes this result when $R$ is a local ring and $S(\frak p)$ satisfies the condition $\frak pR_{\frak p}$ for each $\frak p\in V(\frak a)$. Indeed we deduce that if $\dim R/\frak a=1$, then  $M$ is $\cN\cS$-$\frak a$-cofinite if and only if $\Hom_R(R/\frak a,M)$ and $\Ext_R^1(R/\frak a,M)\in\cN\cS$. Moreover, if $R$ is a local ring of dimension $2$ such that $\cS(\frak p)$ satisfies the condition $C_{\frak pR_{\frak p}}$ for every prime ideal $\frak p$ of $R$ with $\dim R/\frak p\leq 1$, then $M$ is $\cN\cS$-$\frak a$-cofinite if and only if $\Hom_R(R/\frak a,M)$ and $\Ext_R^1(R/\frak a,M)\in\cN\cS$. 

Bahmanpour et all [BNS, Theorem 3.5] showed that if $R$ is a local ring such that $\dim R/\frak a=2$, then $M$ is $\frak a$-cofinite if and only if $\Ext_R^i(R/\frak a,M)$ are finitely generated for $i=0,1,2$. As another conclusion, we generalize this result when $S(\frak p)$ satisfies the condition $\frak pR_{\frak p}$ for each $\frak p\in V(\frak a)$. To be more prcise, we deduce that if $\dim R/\frak a=2$, then $M$ is $\cN\cS$-$\frak a$-cofinite if and only if $\Ext_R^i(R/\frak a,M)\in\cN\cS$ for $i=0,1,2$. Moreover, if $R$ is a local ring of dimension $3$ such that $\cS(\frak p)$ satisfies the condition $\frak pR_{\frak p}$ for every prime ideal $\frak p$ with $\dim R/\frak p\leq 2$, then $M$ is $\cN\cS$-$\frak a$-cofinite if and only if $\Ext_R^i(R/\frak a,M)\in\cN\cS$ for $i=0,1,2$. We prove the following result about local cohomology which generalizes [NS, Theorem 3.7]. For the basic properties and unexplained terminology of local cohomology, we refer the reader to the textbook by Brodmann and Sharp [BS]. 

\begin{Theorem}
Let $(R,\frak m)$ be a local ring, let $\frak a$ be an ideal of $R$ such that $\dim R/\frak a=2$ and let
$S(\frak p)$ satisfy the condition $\frak pR_{\frak p}$ for each $\frak p\in V(\frak a)$. If $n$ is a non-negative integer such that $\Ext_R^i(R/\frak a,M)\in\cN\cS$ for all $i\leq n+1$, then the following conditions are equivalent.
 
{\rm (i)} $H_{\frak a}^i(M)$ is $\cN\cS$-$\frak a$-cofinite for all $i<n$.

{\rm (ii)} $\Hom_R(R/\frak a,H_{\frak a}^i(M))\in\cN\cS$ for all $i\leq n$. 
\end{Theorem}
 
In Section 3, we define a new dimension of modules which is an upper bound of the dimension mentioned in Section 3. For every $R$-module $M$, we define $\overline{\dim}M=\dim R/\Ann_RM$. For every non-negative integer $n$, we denote by $\cD_{\leq n}$, the subcategory of all $R$-modules $M$ satisfying $\overline{\dim} M\leq n$. It is proved that $\cD_{\leq n}$ is Serre for each $n\geq 0$. We show that the subcategory $\frak a$-cofinite modules in $\cD_{\leq 1}$ is Serre, the subcategory of $\frak a$-cofinite modules in $\cD_{\leq 2}$ is abelian. Finally we show that the kernel and cokernel of a homomorphism $f:M\To N$ of $\frak a$-cofinite modules in $\cD_{\leq 3}$ is $\frak a$-cofinite if and only if $(0:_{\Coker f}\frak a)$ is finitely generated.

\section{Cofiniteness with respect to Extension subcategories}

Let $\cC$ be an abelian category and $\cS$ be a subcategory of $\cC$. We denote by $\cS_{\rm sub}$ ($\cS_{\rm quot}$), the smallest subcategory of $\cC$ containing $\cS$ which is closed under subobjects (quotients). These subcategories can be specified as follows:    

$$(\cS)_{\rm sub}=\{N\in\cC|\hspace{0.1cm} N \hspace{0.1cm}{\rm is\hspace{0.1cm} a\hspace{0.1cm} subobject\hspace{0.1cm}
object\hspace{0.1cm} of\hspace{0.1cm} an \hspace{0.1cm}object\hspace{0.1cm} of}\hspace{0.1cm} \cS\};$$

$$(\cS)_{\rm quot}=\{M\in\cC|\hspace{0.1cm} M \hspace{0.1cm}{\rm is\hspace{0.1cm} a\hspace{0.1cm} quotient\hspace{0.1cm}
object\hspace{0.1cm} of\hspace{0.1cm} an \hspace{0.1cm}object\hspace{0.1cm} of}\hspace{0.1cm} \cS.$$

Let $\cT$ be another subcategory of $\cC$. We denote by $\cS\cT$, the extension subcategory of $\cS$ and $\cT$ which is:

$$\cS\cT=\{M\in\cC|\hspace{0.1cm} {\rm there \hspace{0.1cm} exists \hspace{0.1cm} an\hspace{0.1cm}exact\hspace{0.1cm}
 sequence}$$$$
0\To L\To M\To N\To 0\hspace{0.1cm} {\rm with}\hspace{0.1cm} L\in\cS
\hspace{0.1cm}{\rm and}\hspace{0.1cm} N\in\cT\}.$$

For any $n\in \mathbb{N}_0$, we set $\cS^0=\{0\}$ and
$\cS^n=\cS\cS^{n-1}$. In the case where $\cS^2=\cS$, we say
that $\cS$ is {\it closed under extension}. We also define
$(\cS)_{\rm ext}=\bigcup_{n\geq 0}\cS^n$ the smallest
subcategory of $\cC$ containing $\cS$ which is closed under
extension.

A full subcategory $\cS$ of $\cC$ is called {\it Serre} if it is
closed under taking subobjects, quotients and extensions. For any subcategory $\cX$ of $\cC$, we denote by $(\cX)_{\rm Serre}$, the smallest Serre subcategory of $\cA$ containing $\cX$.
\medskip
\begin{Lemma}\label{kan}
Let $\cS$ be a subcategory of $\cC$. Then we have the following conditions.

${\rm (i)}$ $(\cS)_{\rm Serre}=(((\cS)_{\rm sub})_{\rm quot})_{\rm ext}$.

${\rm (ii)}$ $(\cS)_{\rm ext}=\{N\in\cC|\hspace{0.1cm}\textnormal{there exists a finite filtration}\hspace{0.1cm} 0\subset N_1\subset N_2\subset\dots \subset N_n=N\\ \hspace{0.1cm}\textnormal{of subobjects of }\hspace{0.1cm} N\hspace{0.1cm}\textnormal{such that}\hspace{0.1cm} N_i/N_{i-1}\in\cS\hspace{0.1cm}\textnormal{for each}\hspace{0.1cm} 1\leq i\leq n\}$
\end{Lemma}
\begin{proof}
(i) See [K, Proposition 2.4]. (ii) Given $N\in(\cS)_{\rm ext}$, there exists a positive integer $n$ such that $N\in\cS^n$. We prove by induction on $n$ that there exists a filtration 
 $0\subset N_1\subset N_2\subset\dots \subset N_n=N$ of subobjects of $N$ of length $n$ such that $N_i/N_{i-1}\in\cS$ for each $1\leq i\leq n$. If $n=1$, there is nothing to prove. Since $N\in\cS^n$, there exists an exact sequence $0\To N_1\To N\To N/N_1\To 0$ such that $N_1\in\cS$ and $N/N_1\in\cS^{n-1}$. By the induction hypothesis, there exists a finite filtration  $0\subset N_2/N_1\subset \dots \subset N_n/N_1=N/N_1$ such that $N_i/N_{i-1}\in\cS$ for each $2\leq i\leq n$. Consequently, $0\subset N_1\subset N_2\subset\dots\subset N_n=N$ is the desired filtration.  
\end{proof}

\medskip

  In the rest of this paper, $\frak a$ is an ideal of $R$, $\cS$ is a Serre subcategory of $R$-modules and $\cN$ is the subcategory of finitely generated modules in its corresponding module category (such as in the category of $R$-modules or in the category of $R_{\frak p}$-modules for some prime ideal $\frak p$ of $R$). We have the following lemma.
	
	\begin{Lemma}\label{yosh}
	The extension subcategory $\cN\cS$ is Serre for every Serre subcategory $\cS$ of $R$-modules. Moreover $\cS\cN\subseteq \cN\cS$.
	\end{Lemma}	
	\begin{proof}
	See [Y, Corollary 3.3 and Theorem 3.2]. 
	\end{proof}		
	
	An $R$-module $M$ is said to be $\cS$-$\frak a$-{\it cofinite} if $\Supp M\subseteq V(\frak a)$ and $\Ext_R^i(R/\frak a, M)\in\cS$ for all integers $i\geq 0$. An $\cN$-$\frak a$-cofinite module is called $\frak a$-{\it cofinite}. We recall that $\cS$ satisfies the condition $C_{\frak a}$ if for every $R$-module $M$, the following implication holds.
\begin{center}
$C_{\frak a}$: If $\Gamma_{\frak a}(M)=M$ and $(0:_M{\frak a})$ is in $\cS$, then
$M$ is in $\cS$.
\end{center}
	
	\medskip
	\begin{Lemma}\label{sr}
	The Serre subcategory $\cS$ of $R$-modules satisfies the condition $C_{\frak a}$ if and only if it satisfies the condition $C{\sqrt{\frak a}}$.  
	\end{Lemma}
	\begin{proof}
	See [SR, Proposition 2.4].
	\end{proof}
	\medskip
	\begin{Lemma}\label{ext}
Let $N$ be a finitely generated $R$-module and $M$ be an arbitrary $R$-module such that for a non-negative integer $n$, we have $\Ext_R^i(N,M)\in\cS$ for each $i\leq n$. Then $\Ext_R^i(L,M)\in\cS$ for each finitely generated $R$-module $L$ with $\Supp_RL\subseteq\Supp_RN$ and each $i\leq n$. 
\end{Lemma}
\begin{proof}
See [AS, Lemma 2.1]
\end{proof}

	In this section for every $R$-module $M$, $\dim M$ means dimension of $\Supp_RM$ which is the length of the longest chain of prime ideals in $\Supp_RM$.
	\medskip
	
\begin{Proposition}\label{d1}
Let $M$ be an $R$-module of dimension $d$ such that $\Ext_R^i(R/\frak a,M)\in\cN\cS$ for all $i\leq d$ and ${\rm Max} M\subseteq \Supp\cS$ (e.g. if $R$ is a local ring). Then 

${\rm (i)}$ $H_{\frak b}^i(M)$ is $\cN\cS$-$\frak b$-cofinite for all $i\geq 0$ and all ideals $\frak b$ containing $\frak a$ with $\dim R/\frak b\leq 1$ such that $\cS$ satisfies the condition $C_{\frak b}$. More precisely, $H_{\frak b}^i(M)$ is $\cN\cS$-$\frak b$-cofinite for all $i<d$ and $H_{\frak b}^d(M)\in\cS$. 

${\rm (ii)}$ $\Ext_R^i(N,M)\in\cN\cS$ for all $i\geq 0$ and all finitely generated $R$-modules $N$ with $\Supp_R N\subseteq V(\frak a)$ and $\dim N\leq 1$ such that $\cS$ and satisfying the condition $C_{\Ann_R N}$.
\end{Proposition}
\begin{proof}
(i) Assume that $\frak b$ an ideal of $R$ containing $\frak a$ with $\dim R/\frak b\leq 1$. It follows from \cref{ext} that $\Ext_R^i(R/\frak b,M)\in\cN\cS$ for all $i\leq d$. If $d=0$, there exists an exact sequence $$0\To N\To \Hom_R(R/\frak b,M)\To S\To 0 $$ such that $N$ has finite length and $S\in\cS$. Since ${\rm Max} M\subseteq \Supp\cS$, we deduce that $N\in\cS$ so that $\Hom_R(R/\frak b,M)\in\cS$ and since $\cS$ satisfies the condition $C_{\frak b}$, we have $M\in\cS$. This implies that $\Gamma_{\frak b}(M)\in\cS$ and hence $\Gamma_{\frak b}(M)$ is $\cN\cS$-$\frak a$-cofinite (we observe that $H_{\frak b}^i(M)=0$ for all $i>0$). Assume that $d>0$ and we proceed by induction on $i$. If $\dim R/\frak b=0$, since $\Ext_R^i(R/\frak b,M)\in\cN\cS$, using a similar argument as mentioned above, $\Ext_R^i(R/\frak b,M)\in\cS$ for all $i\leq d$; and hence it follows from [AM, Theorem 2.9] that $H_{\frak b}^i(M)\in\cS$ for all $i\geq 0$; and hence the assertion is clear in this case. If $\dim R/\frak b=1$, it follows from [AS, Theorem 3.5] that $H_{\frak b}^i(M)$ is $\cN\cS$-$\frak b$-cofinite for all $i<d$ and $\Hom_R(R/\frak a, H_{\frak b}^d(M))\in\cN\cS$. Since $\Supp_R H_{\frak b}^d(M)\subseteq {\rm Max}R$, using a similar argument as mentioned in the case $d=0$, we deduce that $\Hom_R(R/\frak a, H_{\frak b}^d(M))\in\cS$ and since $\cS$ satisfies the condition $C_{\frak b}$, we deduce that $H_{\frak b}^d(M)\in\cS$. (ii) It follows from \cref{ext} that $\Ext_R^i(R/\Ann_R N,M)\in\cN\cS$ for all $i\geq 0$ if and only if $\Ext_R^i(R/\sqrt{\Ann_R N},M)\in\cN\cS$ for all $i\geq 0$; ad hence we may assume that $\Ann_R N=\sqrt{\Ann_R N}$. Since $\Supp_R N\subseteq V(\frak a)$, we have $\frak a\subseteq \Ann N$; and hence (i) implies that $H_{\Ann_R N}^i(M)$ is $\cN\cS$-$\Ann_R N$-cofinite for all $i\geq 0$. Then it follows from [AS, Theorem 3.5] that $\Ext_R^i(R/\Ann_R N,M)\in\cN\cS$ for all $i\geq 0$. It thus follows from \cref{yosh} and \cref{ext} that $\Ext_R^i(N,M)\in\cN\cS$ for all $i\geq 0$.
\end{proof}

\medskip

For every Serre subcategory $\cS$ of $R$-modules and every $\frak p\in\Spec R$, we denote by $\cS({\frak p})$ the smallest Serre subcategory of $R_{\frak p}$-modules containing $\cS_{\frak p}=\{M_{\frak p}|\hspace{0.1cm} M\in\cS\}$. We have the following lemma.

\begin{Lemma}\label{sqext}
Let $\frak p$ be an ideal of $R$. Then $\cS_{\frak p}$ is closed under subobjects and quotients. In particular, $\cS(\frak p)=(\cS_{\frak p})_{\rm ext}$.
\end{Lemma}
\begin{proof}
Given $M_{\frak p}\in\cS_{\frak p}$ and a submodule $K$ of $M_{\frak p}$, there exists a submodule $N$ of $M$ such that 
$K=N_{\frak p}$. Since $M\in\cS$, we have $N, M/N\in\cS$ and $M_{\frak p}/K=(M/N)_{\frak p}$ which yields $K,M_{\frak p}/K\in\cS_{\frak p}$. The second assertion follows from \cref{kan}.
\end{proof}

By virtue of \cref{sqext}, if $(R,\frak m)$ is a local ring, then we have $\cS(\frak m)=\cS$.

\medskip
\begin{Lemma}\label{ess}
Let $\frak p$ be a prime ideal of $R$ and $X\in\cS_{\frak p}$. Then there exists an $R$-module $S\in\cS$ such that $S$ is an essential $R$-submodule of $S_{\frak p}$ and $X=S_{\frak p}$.
\end{Lemma}
\begin{proof}
Since $X\in\cS_{\frak p}$, there exists $T\in\cS$ such that $X=T_{\frak p}$. If $\phi: T\To T_{\frak p}$ is the canonical homomorphism, then $S=T/\ker \phi$ is the desired module. 
\end{proof}

\medskip 
\begin{Theorem}\label{pdim1}
Let $(R,\frak m)$ be a local ring, let $S(\frak p)$ satisfy the condition $\frak pR_{\frak p}$ for each $\frak p\in V(\frak a)$. If $M$ is an $R$-module of dimension $d$ such that $\Ext_R^i(R/\frak a,M)\in\cN\cS$ for each $i\leq d$, then $\Ext_R^i(N,M)\in\cN\cS$ for each $i\geq 0$ and each finitely generated $R$-module $N$ with $\dim N\leq 1$ and $\Supp_RN\subseteq V(\frak a)$.
\end{Theorem}
\begin{proof}
If $\dim N=0$, then $\Ann_RN$ is $\frak m$-primary. It follows from the assumption and \cref{sr} that $N$ satisfies the condition $C_{\Ann_RN}$; and hence \cref{d1} implies that $\Ext_R^i(N,M)\in\cN\cS$ for each $i\geq 0$. Now, assume that $\dim N=1$,  $\Min\Ann_RN=\{\frak p_1,\dots,\frak p_k\}$ and $\Ext_R^i(R/\frak p_j,M)\in\cN\cS$ for each $i\geq 0$ and each $1\leq j\leq k$. If we set $T=\oplus_{j=1}^kR/\frak p_j$, then $\Ext_R^i(T,M)\in\cN\cS$ for each $i\geq 0$ and we have $\Supp_RT=\Supp_RN$. It thus follows from \cref{ext} that $\Ext_R^i(N,M)\in\cN\cS$ for each $i\geq 0$; and hence we may assume that $N=R/\frak p$ for some $\frak p\in V(\frak a)$ with $\dim R/\frak p=1$.
It is clear that $\Ext_{R_{\frak p}}^i(R_{\frak p}/\frak a R_{\frak p},M_{\frak p})\in\cN\cS({\frak p})$ for each $i\leq d$; and hence it follows from \cref{d1} that $\Ext_{R_{\frak p}}^i(R_{\frak p}/\frak pR_{\frak p},M_{\frak p})\in\cN\cS(\frak p)$ for each $i\geq 0$. Assume that $i\geq 0$ and $L=\Ext_R^i(R/\frak p,M)$. We show that $L\in\cN\cS$. In view of the previous argument there exists a finitely generated submodule $K$ of $L$ such that $L_{\frak p}/K_{\frak p}\in\cS(\frak p)$. Consider the canonical morphism $\phi:L/K\To (L/K)_{\frak p}$ with $\Ker\phi=L_1/K$  where $L_1$ is a submodule of $L$ containing $K$. Clearly $(L_1/K)_{\frak p}=0$ and since $\frak p\subseteq \Ann_R L_1/K$ and $\dim R/\frak p=1$, we have $\Supp_RL_1/K\subseteq V(\frak m)$. Considering $x\in\frak m\setminus \frak p$, there is an exact sequence $$0\To R/\frak p\stackrel{x.}\To R/\frak p\To R/\frak p+xR\To 0$$ which yields the following exact sequence 
$$\Ext_R^i(R/\frak p +xR, M)\To \Ext_R^i(R/\frak p,M)\stackrel{x.}\To\Ext_R^i(R/\frak p,M).$$ We observe that $\frak p +xR$ is $\frak m$-primary and according to the assumption and \cref{sr}, $\cS$ satisfies the condition $C_{\frak p+xR}$, and hence \cref{d1} implies that $\Ext_R^i(R/\frak p +xR,M)\in\cN\cS$. Thus $(0:_Lx)\in\cN\cS$.
Since $(L/K)_{\frak p}\in\cS(\frak p)$, using \cref{kan}, there exists a positive integer $t$ such that $(L/K)_{\frak p}=\cS^t_{\frak p}$.  Without loss of generality, we may assume that $t=2$, the other cases are similar. Then there exists an exact sequence of $R_{\frak p}$-modules $$0\To S'_{\frak p}\To (L/K)_{\frak p}\To S''_{\frak p}\To 0$$ such that $S',S''\in\cS$. Using \cref{ess} we may assume that $S'$ and $S''$ are $R$-submodule of $S'_{\frak p}$ and $S''_{\frak p}$, respectively. Since $L/L_1$ is an essential $R$-submodule of $(L/L_1)_{\frak p}=(L/K)_{\frak p}$, the $R$-submodule $S'\cap L/L_1$ is nonzero and $(S'\cap L/L_1)_{\frak p}=S'_{\frak p}$ and replacing $S'$ by $S'\cap L/L_1$, we may assume that $S'$ is a submodule of $L/L_1$. Suppose that $S'=L_2/L_1$ for submodule $L_2$ of $L$ containing $L_1$ and we have $(L/L_2)_{\frak p}=S''_{\frak p}$. Since $(L_1/K)_{\frak p}=0$, we have $\Supp_RL_1/K\subseteq V(\frak m)$ and applying the functor $\Hom_R(R/xR,-)$ to the exact sequence $0\To K\To L_1\To L_1/K\To 0$ and the previous argument, we deduce that $\Hom_R(R/xR,L_1/K)\in\cN\cS$. Then there is an exact sequence $$0\To F\To\Hom_R(R/xR,L_1/K)\To S\To 0$$ such that $F$ is finitely generated and $S\in\cS$. Since $\Supp_R\Hom_R(R/xR, L_1/K)\subseteq V(\frak m)$, the module $F$ has finite length and hence $F\in\cS$ as $R/\frak m\in\cS$. This implies that $\Hom_R(R/xR,L_1/K)\in\cS$ and hence $\Hom_R(R/\frak m,L_1/K)\in\cS$. Since $\cS$ satisfies the condition $C_{\frak m}$, we deduce that $L_1/K\in\cS$ so that $L_1\in\cN\cS$ and consequently $L_2\in\cN\cS$. Consider the canonical homomorphism $\phi_1:L/L_2\To (L/L_2)_{\frak p}$ with $L_3/L_2=\Ker \phi_1$. Thus $(L_3/L_2)_{\frak p}=0$ so that $\Supp_RL_3/L_2\subseteq V(\frak m)$. From the induced essential monomorphism $L/L_3\To (L/L_2)_{\frak p}=(L/L_3)_{\frak p}$ and using a similar argument as mentioned above, we may assume that $S''$ is a submodule of $L/L_3$ and so there exists a submodule $L_4$ of $L$ containing $L_3$ such that $S''=L_4/L_3$ and $(L/L_4)_{\frak p}=0$. Therefore we have $\Supp_RL/L_4\subseteq V(\frak m)$.  The exact sequence $0\To L_2\To L_3\To L_3/L_2\To 0$, \cref{yosh} and the fact that $L_2\in\cN\cS$ imply that $\Hom_R(R/xR,L_3/L_2)\in\cN\cS$ and so using a similar argument as mentioned above, we deduce that $L_3/L_2\in\cS$ so that $L_3\in\cN\cS$ and hence $L_4\in\cN\cS$. Now applying $\Hom_R(R/xR,-)$ to the exact sequence $0\To L_4\To L\To L/L_4\To 0$ and using again a similar argument as mentioned before, we deduce that $L\in\cN\cS$.     
\end{proof}
\medskip
The following corollary generalizes a result due to Melkersson [M2, Theorem 2.3].
\begin{Corollary}\label{t22}
Let $(R,\frak m)$ be a local ring, let $\dim R/\frak a\leq 1$ and let
$S(\frak p)$ satisfy the condition $\frak pR_{\frak p}$ for each $\frak p\in V(\frak a)$. If $M$ is an $R$-module such that $\Supp_RM\subseteq V(\frak a)$ and $\Hom_R(R/\frak a,M), \Ext_R^1(R/\frak a,M)\in\cN\cS$, then $M$ is $\cN\cS$-$\frak a$-cofinite.
\end{Corollary}
\begin{proof}
The assertion follows immediately from \cref{pdim1}.
\end{proof}

  \medskip
\begin{Theorem}\label{d2}
Let $(R,\frak m)$ be a local ring, let $S(\frak p)$ satisfy the condition $\frak pR_{\frak p}$ for each $\frak p\in V(\frak a)$. If $M$ is an $R$-module of dimension $d$ such that $\Ext_R^i(R/\frak a,M)\in\cN\cS$ for each $i\leq d$, then $\Ext_R^i(N,M)\in\cN\cS$ for each $i\geq 0$ and each finitely generated $R$-module $N$ with $\dim N\leq 2$ and $\Supp_RN\subseteq V(\frak a)$.
\end{Theorem}
\begin{proof}
Similar to the proof of \cref{pdim1}, we may assume that $N=R/\frak p$ for some $\frak p\in V(\frak a)$ with $\dim R/\frak p=2$ and $\Ext_{R_{\frak p}}^i(R_{\frak p}/\frak pR_{\frak p},M_{\frak p})\in\cN\cS(\frak p)$ for each $i\geq 0$. Assume that $i\geq 0$ and $L=\Ext_R^i(R/\frak p,M)$ and we show that $L\in\cN\cS$. There exists a finitely generated submodule $K$ of $L$ such that $L_{\frak p}/K_{\frak p}\in\cS(\frak p)$. Consider the canonical morphism $\phi_{L/K}:L/K\To (L/K)_{\frak p}$ with $\Ker\phi_{L/K}=L_1/K$ where $L_1$ is a submodule of $L$. Clearly $(L_1/K)_{\frak p}=0$ and since $\frak p\subseteq \Ann_R L_1/K$ and $\dim R/\frak p=2$, every $\frak q\in\Supp_RL_1/K$ with $\dim R/\frak q=1$ is a minimal prime ideal of $\Ann_R L_1/K$. Then the set $$\cT=\{\frak q\in\Supp_RL_1/K|\hspace{0.1cm} \dim R/\frak q=1\}$$ is finite. Assume that $\cT=\{\frak q_1,\dots,\frak q_n\}$. Then $\frak p\subsetneq\cap_{j=1}^n\frak q_j$ and so there $x\in\cap_{i=1}^n\frak q_i\setminus \frak p$. Then there is an exact sequence $$0\To R/\frak p\stackrel{x.}\To R/\frak p\To R/\frak p+xR\To 0$$ which yields the following exact sequence 
$$\Ext_R^i(R/\frak p +xR, M)\To \Ext_R^i(R/\frak p, M)\stackrel{x.}\To\Ext_R^i(R/\frak p, M).$$ Since $\dim R/\frak p+xR\leq 1$, it follows from \cref{pdim1} that $\Ext_R^i(R/\frak p+xR,M)\in\cN\cS$; and hence $(0:_Lx)\in\cN\cS$. Now assume that $L_2/K=\Ker\phi_{L_1/K}$ where $\phi_{L_1/K}:L_1/K\To (L_1/K)_{\frak q_1}$. Then we have $(L_2/K)_{\frak q_1}=0$. Continuing this way, assume that $L_{i+1}/K=\Ker\phi_{L_i/K}$ where $\phi_{L_i/K}:L_i/K\To (L_i/K)_{\frak q_i}$ for $1\leq i\leq n$, we deduce that $\Supp_RL_{i}/K\subseteq \{\frak q_i,\dots \frak q_n,\frak m\}$ for each $1\leq i\leq n$, $\Supp_RL_{n+1}/K\subseteq \{\frak m\}$ and we have a chain of submodules $$L_{n+1}\subseteq L_n\subseteq \dots\subseteq L_1\subseteq L=L_0$$ such that the induced morphism $\overline{\phi_{L_i/K}}:L_i/ L_{i+1}\To (L_i/K)_{\frak q_i}$ is an essential monomorphism for each $0\leq i\leq n$. We observe that $\Supp_RL_i/L_{i+1}\subseteq \{\frak q_i,\dots \frak q_n,\frak m\}$ for each $1\leq i\leq n$. Since $(0:_{L}x)\in\cN\cS$, we have $(0:_{L_{n+1}}x)\in\cN\cS$. Now applying the functor $\Hom_R(R/xR,-)$ to the exact sequence $0\To K\To L_{n+1}\to L_{n+1}/K\To 0$, we deduce that $\Hom_R(R/xR, L_{n+1}/K)\in\cN\cS$ and so there is an exact sequence of $R$-modules  $$0\To F\To\Hom_R(R/xR, L_{n+1}/K)\To S\To 0$$ such that $F$ is finitely generated and $S\in\cS$. Using the same argument as mentioned in the proof of \cref{pdim1}, we deduce that $L_{n+1}/K\in\cS$ so that $L_{n+1}\in\cN\cS$. Since $\Hom_R(R/xR,L_n)\in\cN\cS$, $\Hom_R(R/xR,L_n)_{\frak q_n}\in\cN\cS(\frak q_n)$ so that $\Hom_R(R/xR,L_n/K)_{\frak q_n}\in\cN\cS(q)$ and a similar proof as mentioned above gives $(L_n/K)_{\frak q_n}=(L_n/L_{n+1})_{\frak q_n}\in\cS(\frak \frak q_n)$. To be convenient, set $\frak q=\frak q_n$. in view of \cref{sqext}, we may assume that $(L_n/L_{n+1})_{\frak q}=\cS_{\frak q}^2$ and the other cases are similar. Then there exists an exact sequence 
$$0\To S'_{\frak q}\To (L_n/L_{n+1})_{\frak q}\To S''_{\frak q}\To .$$ Using \cref{ess}, we may assume that $S'$ and $S''$ are $R$- submodules of $S'_{\frak p}$ and $S''_{\frak p}$, respectively. Since $L_n/L_{n+1}$ is an essential submodule of $(L_n/L_{n+1})_{\frak q}$, the submodule $S'\cap L_n/L_{n+1}$ is nonzero and $(S'\cap L_n/L_{n+1})_{\frak p}=S'_{\frak p}$; and hence replacing $S'$ by $S'\cap L_n/L_{n+1}$, we may assume that $S'$ is a submodule of $L_n/L_{n+1}$. Assume that $S'=L'_n/L_{n+1}$ for some submodule $L'_n$ of $L_n$ and so $L'_n\in\cN\cS$. Thus $(L_n/L_n')_{\frak q}=S''_{\frak q}$. Consider the canonical exact sequence $0\To L''_n/L'_n\To L_n/L'_n\To (L_n/L'_n)_{\frak q}$ which forces $\Supp_R L''_n/L'_n\subseteq V(\frak m)$. Since $\Hom_R(R/xR,L''_n)\in\cN\cS$ and $L'_n\in\cN\cS$, using a similar argument as above, we deduce that $L''_n\in\cN\cS$. The essential monomorphism $L_n/L''_n\To (L_n/L'_n)_{\frak q}$ and a similar argument as above imply that $S''$ is a submodule of $L_n/L''_n$ and so $S''=L''/L''_n$ for some submodule $L''$ of $L_n$. The fact that $L''_n\in\cN\cS$ implies that $L''\in\cN\cS$. Since $(L_n/L'')_{\frak q}=0$, we deduce that $\Supp_RL_n/L''\subseteq V(\frak m)$. A similar argument as above implies that $L_n/L''\in\cS$, and hence $L_n\in\cN\cS$. Continuing this manner we deduce that $L=L_1\in\cN\cS$. Therefore, we may assume that $\phi_{L/K}$ is an essential monomorphism. According to \cref{sqext}, there exists a positive integer $t$ such that $(L/K)_{\frak p}\in(\cS_{\frak p})^t$. Set $A=L/K$. Without loss of generality, we may assume that $t=2$ and so by a similar argument as mentioned above, there exists an exact sequence of $R_{\frak p}$-modules $0\To S'_{\frak p}\To A_{\frak p}\To S''_{\frak p}\To 0$ such $S'$ is a submodule of $L/K$ and $S''$ is an $R$-submodule of $A/S'$ and so $((A/S')/S'')_{\frak p}=0$; and hence $\frak p\notin\Supp_R(A/S')/S''$. Since $\frak p\subseteq \Ann_R (A/S')/S''$ and $\dim R/\frak p=2$, every $\frak q\in\Supp_R(A/S')/S''$ with $\dim R/\frak q=1$ is a minimal prime ideal of $\Ann_R (A/S')/S''$; and hence the set $$\cU=\{\frak P\in\Supp_R(A/S')/S''|\hspace{0.1cm} \dim R/\frak P=1\}$$ is finite. Assume that $\cU=\{\frak p_1,\dots\frak p_m\}$ and set $D=(A/S')/S''$. We notice that $\Supp_RD=\{\frak p_1,\dots,\frak p_m,\frak m\}$. It follows from \cref{d1} that $D_{\frak p_j}\in\cN\cS(\frak p_j)$ for each $1\leq j\leq m$. Then for each $j$, there exists an exact sequence $$0\To N_{j_{\frak p_j}}\To D_{\frak p_j}\To (D/N)_{\frak p_j}\To 0$$ such that $N_j$ is a finitely generated $R$-module and $(D/N)_{\frak p_j}\in (S_{\frak p_j})_{\rm ext}$. In view of the preceding arguments, we may assume that $S\in\cS_{\frak p_j}$ so that there exists $S_j\in\cS$ such that $S=S_{j_{\frak p_j}}$. Furthermore, we may assume that $S_j$ is a submodule of $D/N_j$ and $((D/N_j)/S_j)_{\frak p_j}=0$. Then there exists a submodule $X_j$ of $D$ such thst $S_j=X_j/N_j$ and so this implies that $X_j\in\cN\cS$ for each $1\leq j\leq m$. Now, set $X=X_1+\dots+ X_m$ and so clearly $X\in\cN\cS$ and $\Supp D/X\subseteq V(\frak m)$. We notice that there exists a submodule $A'$ of $A$ such that $S''=A'/S'$  and so $A'\in\cS$. On the other hand, there exists a submodule $L'$ of $L$ such that $A'=L'/K$ which implies that $L'\in\cN\cS$ and $D=L/L'$. Moreover, since $X$ is a submodule of $D$, there exists a submodule $L_1$ of $L$ containing $L'$ such that $X=L_1/L'$ and hence $D/X=L/L_1$. We observe by \cref{yosh} that $L_1\in\cN\cS\cN\cS\subseteq \cN\cS$. Considering the following exact sequence $$0\To L_1\To L\To D/X\to 0,$$ we have $\Hom_R(R/xR,D/X)\in\cN\cS$ and using a similar proof as mentioned before, we deduce that $L/L_1=D/X\in\cS$. Now the fact that $L_1\in\cN\cS$ forces $L\in\cN\cS\cS\subseteq\cN\cS$.  
\end{proof}

The following corollary generalizes [BNS, Theorem 3.5].
\medskip

\begin{Corollary}\label{t21}
Let $(R,\frak m)$ be a local ring, let $\dim R/\frak a\leq 2$ and let
$S(\frak p)$ satisfy the condition $\frak pR_{\frak p}$ for each $\frak p\in V(\frak a)$. If $M$ is an $R$-module such that $\Supp_RM\subseteq V(\frak a)$ and $\Ext_R^i(R/\frak a,M)\in\cN\cS$ for $i=0,1,2$, then $M$ is $\cN\cS$-$\frak a$-cofinite.
\end{Corollary}
\begin{proof}
The assertion follows immediately from \cref{d2}. 
\end{proof}

The subcategory of all $R$-modules of finite support is denoted by $\cF$. It is clear that $\cF$ is Serre and it satisfies the condtion $C_{\frak a}$ for every ideal $\frak a$ of $R$.
\medskip
\begin{Corollary}\label{f}
Let $(R,\frak m)$ be a local ring, let $\dim R/\frak a=2$ and let $M$ be an $R$-module such that $\Supp_RM\subseteq V(\frak a)$. If $\Ext_R^i(R/\frak a,M)\in\cN\cF$ for $i=0,1,2$, then $M$ is $\cN\cF$-$\frak a$-cofinite.
\end{Corollary}
\begin{proof}
It is clear that $\cF(\frak p)$ satisfies the condition $C_{\frak pR_{\frak p}}$ for every prime ideal $\frak p$ of $R$. Therefore, the assertion follows immediately from \cref{d2}. 
\end{proof}
\medskip

\begin{Theorem}\label{loc}
Let $(R,\frak m)$ be a local ring, let $\dim R/\frak a=2$ and let
$S(\frak p)$ satisfy the condition $\frak pR_{\frak p}$ for each $\frak p\in V(\frak a)$. If $n$ is a non-negative integer such that $\Ext_R^i(R/\frak a,M)\in\cN\cS$ for all $i\leq n+1$, then the following conditions are equivalent.
 
{\rm (i)} $H_{\frak a}^i(M)$ is $\cN\cS$-$\frak a$-cofinite for all $i<n$.

{\rm (ii)} $\Hom_R(R/\frak a,H_{\frak a}^i(M))\in\cN\cS$ for all $i\leq n$. 
\end{Theorem}
\begin{proof}
We proceed by induction on $n$. The case $n=0$ is clear. If we consider $\overline{M}=M/\Gamma_{\frak a}(M)$, there exists an exact sequence of $R$-modules $0\To \overline{M}\To E\To Q\To 0$ such that $E$ is injective with $\Gamma_{\frak a}(E)=0$. This eaxct sequence implies that $H_{\frak a}^i(Q)\cong H_{\frak a}^{i+1}(M)$ and $\Ext_R^i(R/\frak a,Q)\cong\Ext_R^i(R/\frak a,\overline{M})$ for each $i\geq 0$; furthermore $\Hom_R(R/\frak a,H_{\frak a}^1(M))\cong\Ext_R^1(R/\frak a,\overline{M})$. Moreover for each $i\geq 1$ we have an exact sequence of $R$-modules,
$$\Ext_R^i(R/\frak a,\Gamma_{\frak a}(M))\To\Ext_R^i(R/\frak a,M)\To\Ext_R^i(R/\frak a,\overline{M})$$$$\To\Ext_R^{i+1}(R/\frak a,\Gamma_{\frak a}(M))\To\Ext_R^{i+1}(R/\frak a,M)\hspace{0.3cm}(\dag_i).$$
If $n=1$, the exact sequence $(\dag_1)$ and \cref{t21} imply that (i) and (ii) are equivalent. For $n>1$, since $\Gamma_{\frak a}(M)\in\cN\cS$, it follows from $(\dag_i)$ and the previous isomorphisms that $\Ext_R^j(R/\frak a,Q)\in\cN\cS$ for all $i\leq n$ so that (i) and (ii) are equivalent for $Q$ and non-negative integer $n-1$. Now, using again the previous isomorphisms, the conditions (i) and (ii) are equivalent for $M$ and non-negative integer $n$.
\end{proof}

\medskip
\begin{Corollary}
Let $R$ be an local ring with $\dim R/\frak a\leq 2$, let $n$ be a non-negative integer such that $\Ext_R^i(R/\frak a,M)\in\cN\cF$ for all $i\leq n+1$. Then the following conditions are equivalent.
 
{\rm (i)} $H_{\frak a}^i(M)$ is $\cN\cF$-$\frak a$-cofinite for all $i<n$.

{\rm (ii)} $\Hom_R(R/\frak a,H_{\frak a}^i(M))\in\cN\cF$ for all $i\leq n$. 
\end{Corollary}
\begin{proof}
Since $\cF(\frak p)$ satisfies the condition $C_{\frak pR_{\frak p}}$ for every prime ideal $\frak p$ of $R$, the results is obtained by \cref{loc}.
\end{proof}

\medskip
 Given an arbitrary Serre subcategory $\cS$ of $R$-modules, we say that $R$ admits the condition $P_n^{\cS}(\frak a)$ if for every $R$-module $M$, the following implication holds:
\begin{center}
$P_n^{\cS}(\frak a):\hspace{0.3cm}$ {\it If $\Ext_R^i(R/\frak a, M)\in\cS$ for all $i\leq n$ and $\Supp(M)\subseteq V(\frak a)$,
then $M$ is $\cS$-$\frak a$-cofinite.}
 \end{center}

\medskip

\begin{Theorem}\label{td}
Let $R$ be a ring of dimension $d\geq 1$ admitting the condition $P^{\cS}_{d-1}(\frak a)$ for all ideals $\frak a$ of dimension $d-1$ (i.e. $ \dim R/\frak a\leq d-1$), then $R$ admits $P^{\cS}_{d-1}(\frak a)$ for all ideals $\frak a$ of $R$. 
\end{Theorem}
\begin{proof} 
Let $M$ be an $R$-module and $\frak a$ be an arbitrary ideal of $R$ such that $\Supp_RM\subseteq V(\frak a)$ and $\Ext_R^i(R/\frak a,M)\in\cN\cS$ for all $i\leq d-1$. If there exists some positive integer $n$ such that $\frak a^n=0$, then $M=(0:_M\frak a^n)$. On the other hand, if $\frak a=(a_1,\dots,a_t)$, there is an exact sequence of $R$-modules $0\To (0:_M\frak a)\To (0:_M\frak a^n)\stackrel{f}\To a_1(0:_M\frak a^n)\oplus\dots a_t(0:_M\frak a^n)$ where $f(x)=(a_1x,\dots,a_tx)$ . It is clear that $a_i(0:_M\frak a^n)$ is a submodule of $(0:_M\frak a^{n-1})$ for each $i$ and hence, since $\cN\cS$ is Serre, an easy induction on $n$ implies that $(0:_M\frak a^n)\in\cN\cS$. Since $V(\frak a)=\Spec R$, the module $M$ is $\cN\cS$-cofinite. Now, assume that $\frak a$ is not nilpotent and so there is a positive integer $n$ such that $\Gamma_{\frak a}(R)=(0:_R\frak a^n)$. Taking $\overline{R}=R/\Gamma_{\frak a}(R)$ and $\overline{M}=M/(0:_M\frak a^n)$, it is clear that $\overline{M}$ is an $\overline{R}$-module and since $\Gamma_{\frak a}(\overline{R})=0$, $\frak a$ contains an $\overline{R}$-regular element so that $\dim R/\frak a+\Gamma_{\frak a}(R)\leq d-1$. We observe that $\Supp_R\overline{M}\subseteq V(\frak a+\Gamma_{\frak a}(R))$ and it follows from \cref{ext} that $\Ext_R^i(R/\frak a+\Gamma_{\frak a}(R),\overline{M})\in\cN\cS$ for all $i\leq d-1$. Thus the assumption implies that $\overline{M}$ is $\cN\cS$-$\frak a+\Gamma_{\frak a}(R)$-cofinite. We now show that $\Ext_{\overline{R}}^i(\overline{R}/\frak a\overline{R},\overline{M})\in\cN\cS$ for each $i\geq 0$. Consider the Grothendieck spectral sequence $$E_2^{p,q}:=\Ext_{\overline{R}}^p(\Tor_q^R(\overline{R},R/\frak a+\Gamma_{\frak a}(R)),\overline{M})\Rightarrow H^{p+q}=\Ext_R^{p+q}(R/\frak a+\Gamma_{\frak a}(R),\overline{M}).$$
 For $i=0$, we have $\Hom_{\overline{R}}(\overline{R}/\frak a\overline{R},\overline{M})=\Hom_{\overline{R}}(\overline{R}\otimes_RR/\frak a+\Gamma_{\frak a}(R),\overline{M})\cong \Hom_R(R/\frak a+\Gamma_{\frak a}(R),\overline{M})\in\cN\cS$. Now, assume that $i>0$ and the result has been proved for all values smaller than $i$. Then $E_2^{p,0}=\Ext_{\overline{R}}^p(\overline{R}\otimes_RR/\frak a+\Gamma_{\frak a}(R),\overline{M})\in\cN\cS$ for all $0\leq p<i$.  Since $\Supp_{\overline{R}}\Tor_q^R(\overline{R},R/\frak a+\Gamma_{\frak a}(R))\subseteq V(\frak a\overline{R})$, It follows from \cref{ext} that $E_2^{p,q}\in\cN\cS$ for all $0\leq p<i$ and all $q\geq 0$. The exact sequence $E_2^{i-2,1}\To E_2^{i,0}\To E_3^{i,0}\To 0$ and the induction hypothesis imply that $E_2^{i,0}\in\cN\cS$ if $E_3^{i,0}\in\cN\cS$. Continuing this manner, we deduce that $E_2^{i,0}\in\cN\cS$ if $E_{\infty}^{i,0}\in\cN\cS$. But there are the following filtration 
$$0=\Phi^{i+1}H^i\subset\dots\subset\Phi^1H^i\subset \Phi^0H^i\subset H^i$$ such that $E_{\infty}^{i,0}\cong \Phi^iH^i/\Phi^{i+1}H^i =\Phi^iH^i$ is a submodule of $H^i=\Ext_R^i(R/\frak a+\Gamma_{\frak a}(R),\overline{M})$; and hence it is in $\cN\cS$. Therefore $\Ext_{\overline{R}}^i(\overline{R}/\frak a\overline{R},\overline{M})\in\cN\cS$ for all $i\geq 0$. Now consider the Grothendieck spectral sequence $$E_2^{p,q}:=\Ext_{\overline{R}}^p(\Tor_q^R(\overline{R},R/\frak a),\overline{M})\Rightarrow H^{p+q}=\Ext_R^{p+q}(R/\frak a,\overline{M}).$$  
Using again \cref{ext}, we deduce $E_2^{p,q}\in\cN\cS$ for all $p,q\geq 0$. For any $r>2$, the $\overline{R}$-module $E_r^{p,q}$ is a subquotient of $E_{r-1}^{p,q}$ and so an easy induction yields that $E_r^{p,q}\in\cN\cS$ for all $r\geq 2$  so that $E_{\infty}^{p,q}\in\cN\cS$ for all $p,q\geq 0$ . For any $0\leq t\leq n$, there is a finite filtration 
$$0=\Phi^{t+1}H^t\subset \Phi^tH^t\subset\dots\subset\Phi^1H^t\subset \Phi^0H^t\subset H^t$$ 
such that $\Phi^pH^t/\Phi^{p+1}H^t\cong E_{\infty}^{p,t-p}$ where $0\leq p\leq t$. Since $E_{\infty}^{p,t-p}\in\cN\cS$ for all $0\leq p\leq t$ and $t\geq 0$, we deduce that $H^t\in\cN\cS$ for all $t\geq 0$; and hence $\overline{M}$ is $\cN\cS$-$\frak a$-cofinite. On the other hand, since $(0:_M\frak a^n)\in\cN\cS$, we conclude that $M$ is $\cN\cS$-$\frak a$-cofinite.
\end{proof}
\medskip

\begin{Corollary}
Let $R$ be a local ring of dimension $2$ such that $\cS$ satisfies the condition $C_{\frak a}$ for every ideal $\frak a$ of $R$ with $\dim R/\frak a\leq 1$. Then $R$ admits the condition $P_1^{\cN\cS}(\frak a)$ for every ideal $\frak a$ of $R$. 
\end{Corollary}
\begin{proof}
It follows from [AS, Theorem 3.2] and \cref{t22} that $R$ admits the condition $P_1^{\cN\cS}(\frak a)$ for all ideals with $\dim R/\frak a\leq 1$. Now, the result follows from \cref{td}.
\end{proof}

\medskip
\begin{Corollary}
Let $R$ be a local ring of dimension $2$ such that $\cS(\frak p)$ satisfies the condition $\frak pR_{\frak p}$ for every prime ideal $\frak p$ with $\dim R/\frak p\leq 1$. Then $R$ admits the condition $P_2^{\cN\cS}(\frak a)$ for every ideal $\frak a$ of $R$. 
\end{Corollary}
\begin{proof}
It follows from \cref{t22} that $R$ admits the condition $P_1^{\cN\cS}(\frak a)$ for all ideals $\frak a$ with $\dim R/\frak a\leq 1$. Now, the result follows from \cref{td}.
\end{proof}
\medskip

\begin{Corollary}
Let $R$ be a local ring of dimension $3$ such that $\cS(\frak p)$ satisfies the condition $\frak pR_{\frak p}$ for every prime ideal $\frak p$ with $\dim R/\frak p\leq 2$. Then $R$ admits the condition $P_2^{\cN\cS}(\frak a)$ for every ideal $\frak a$ of $R$. 
\end{Corollary}
\begin{proof}
It follows from \cref{t21} that $R$ admits the condition $P_2^{\cN\cS}(\frak a)$ for all ideals with $\dim R/\frak a\leq 2$. Now, the result follows from \cref{td}.
\end{proof}

\medskip
\begin{Corollary}
Let $R$ be a local ring of dimension $3$. Then $R$ admits the condition $P_2^{\cN\cF}(\frak a)$. 
\end{Corollary}
\begin{proof}
It follows from \cref{f} that $R$ admits the condition $P_2^{\cN\cF}(\frak a)$ for all ideals with $\dim R/\frak a=2$. Now, the result follows from \cref{td}.
\end{proof}

\section{Cofiniteness with respect to a new dimension}

For every $R$-module $M$, it is clear that $\Supp_RM\subseteq V(\Ann_RM)$ and for the case where $M$ is finitely generated, they are equal. We define the {\it upper dimension} of $M$ and we denote it by $\overline{\dim}M$ which is $\overline{\dim}M=\dim R/\Ann_RM$. Clearly $\dim M\leq \overline{\dim}M$. We first recall some results which are needed in this section.  

\begin{Lemma}\label{del}
Let $S$ be a finitely generated $R$-algebra and let $M$ be an $S$-module. Then $M$ is $\frak a$-cofinite if and only if $M$ is $\frak a S$-cofinite (as an $S$-module).
\end{Lemma}
\begin{proof}
See [DM, Proposition 2].
\end{proof}
\medskip
\begin{Lemma}
Let $S$ be a finitely generated $R$-algebra and let $M$ be an $S$-module. Then $M$ satisfies the condition $P^{\cN}_n(\frak a)$ if and only if $M$ satisfies the condition $P_n^{\cN}(\frak aS)$. 
\end{Lemma}
\begin{proof}
See [KS, Proposition 2.15].
\end{proof}

\medskip
\begin{Lemma}\label{mel}
Let $M$ be an $R$-module such that $\Supp_RM\subseteq V(\frak a)$. Then $M$ is Artinian and $\frak a$-cofinite if and only if $(0:_M{\frak a})$ has finite length.
\end{Lemma}
\begin{proof}
See [M1, Proposition 4.1].
\end{proof}
For every non-negative integer $n$, we denote by $\cD_{\leq n}$, the subcategory of all $R$-modules $M$  such that $\overline{\dim} M\leq n$. We also denote by $\cG$, the subcategory of all $R$-modules $F$ such that $V(\Ann_RF)$ is a finite set. 

\medskip
\begin{Lemma}\label{ser}
If $0\To N\To M\To M/N\To 0$ is an exact sequence of $R$-modules, Then $V(\Ann_RM)=V(\Ann_RN)\cup V(\Ann_RM/N)$.
\end{Lemma}
\begin{proof}
Since $\Ann_RM\subseteq \Ann_RN\cap\Ann_R M/N$, we conclude that $V(\Ann_RN)\cup V(\Ann_RM/N)\subseteq V(\Ann_RM)$. Now assume that $\frak p\in V(\Ann_RM)$ and $\frak p\notin V(\Ann_RN)$. Then there exists $r\in\Ann_RN\setminus \frak p$ and so for every $x\in\Ann_R M/N$, we have $rx\in\Ann_RM\subseteq \frak p$ which implies that $x\in\frak p$. Consequently, $\frak p\in V(\Ann_RM/N)$. 
\end{proof}

\medskip

\begin{Corollary}
The following conditions hold.

{\rm (i)} For every non-negative integer $n$, the suncategory $\cD_{\leq n}$ is Serre.

{\rm (ii)} The subcategory $\cG$ is Serre.
\end{Corollary}
\begin{proof}
The proof is is straightforward by \cref{ser}.
\end{proof}

\medskip 

\begin{Proposition}
Let $M$ be an $R$-module with $\overline{\dim}M=2$. Then $M$ is $\frak a$-cofinite if and only if $\Hom_R(R/\frak a,M)$ and $\Ext_R^1(R/\frak a,M)$ are finitely generated.  
\end{Proposition}
\begin{proof}
Let $\overline {R}=R/\Ann_RM$. Using [KS, Proposition 2.15] we may assume that $\dim R=2$. Now, the result follows from [NS, Corollary 2.4].
\end{proof}

\medskip 

\begin{Proposition}
Let $M$ be an $R$-module with $\overline{\dim}M=3$. Then $M$ is $\frak a$-cofinite if and only if $\Ext_R^i(R/\frak a,M)$ is finitely generated for $i=0,1,2$.  
\end{Proposition}
\begin{proof}
Let $\overline {R}=R/\Ann_RM$. Using [KS, Proposition 2.15] we may assume that $\dim R=3$. Now, the result follows from [NS, Corollary 2.5].
\end{proof}
\medskip

\begin{Proposition}\label{dim1}
The subcategory of $\frak a$-cofinite modules in $\cD_{\leq 1}$ is Serre.
\end{Proposition}
\begin{proof}
 If $\overline{\dim} M=0$, then $\dim M=0$ and so the module $(0:_M\frak a)$ has finite length. Thus it follows from \cref{mel} that $M$ is $\frak a$-cofinite. Now, assume that $\overline{\dim} M=1$ and so $\dim \overline{R}=1$ where $\overline{R}=R/\Ann_RM$. It is clear that $(0:_M\frak a\overline{R})=(0:_M\frak a)$ is finitely generated and $\Supp_{\overline{R}}M\subseteq V(\overline{\frak a})$. It follows from [M1, Proposition 4.5] that $M$ is $\overline{\frak a}$-cofinite. Finally, using \cref{del}, $M$ is $\frak a$-cofinite. The second assertion is straightforward. 
\end{proof}

\medskip
\begin{Corollary}
Let $M\in\cN\cG$ with $\Supp_RM\subseteq V(\frak a)$. Then $M$ is $\frak a$-cofinite if and only if $(0:_M\frak a)$ is finitely generated. 
\end{Corollary}
\begin{proof}
Since $M\in\cN\cF$, there exists an exact sequence of $R$-modules $0\To N\To M\To F\To 0$ such that $N$ is finitely generated and $F$ has finite support. We notice that $(0:_F\frak a)$ is finite, and it suffices to show that $F$ is $\frak a$-cofinite and so we may assume that $V(\Ann_RM)$ is a finite set so that $\overline{\dim}M\leq 1$. It follows from [M1, Proposition 4.5] that $M$ is $\overline{\frak a}$-cofinite where $\overline{R}=R/\Ann_RM$ and $\overline{\frak a}=\frak a\overline{R}$. Now \cref{del} implies that $M$ is $\frak a$-cofinite.  
\end{proof}
\medskip

\begin{Proposition}
 The subcategory of $\frak a$-cofinite modules in $\cD_{\leq 2}$ is abelian.
\end{Proposition}
\begin{proof}
Assume that $f:M\To N$ be a morphism of $\frak a$-cofinite modules and assume that $K=\Ker f, I=\Im f$ and $C=\Coker f$. 
 The assumption implies that $(0:_I\frak a)=(0:_I\frak a\overline{R})$ is finitely generated $\overline{R}$-module where $\overline{R}=R/\Ann_RM$. If $\dim \overline{R}=0$, the module $(0:_I\frak a)$ has finite length and so $I$ is Artinian. Now, \cref{mel} implies that $I$ is $\frak a$-cofinite. If $\dim\overline{R}=1$, it follows from \cref{dim1} that $I$ is $\frak a\overline{R}$-cofinite as $M$ is $\frak a\overline{R}$-cofinite; and hence $I$ is $\frak a$-cofinite using \cref{del}. If $\dim\overline{R}=2$, it follows from [NS, Corollary 2.6] that $I$ is $\frak a\overline{R}$-cofinite and so is $\frak a$-cofinite by using \cref{del}. Now, using the exact sequences of $R$-modules $$0\To K\To M\To I\To 0;$$$$0\To I\To N\To C\To 0;$$ it is straightforward to show that $K$ and $C$ are $\frak a$-cofinite modules.
\end{proof}

\medskip

\begin{Proposition}
The kernel and cokernel of a homomorphism $f:M\To N$ of $\frak a$-cofinite modules in $\cD_{\leq 3}$ is $\frak a$-cofinite if and only if $(0:_{\Coker f}\frak a)$ is finitely generated.
\end{Proposition}
\begin{proof}
By the assumption we have $\dim R/\Ann_RM\leq 3$ and also $\dim R/\Ann_RN\leq 3$. If we put $\frak b=\Ann_RM\cap \Ann_RN$ and $\overline{R}=R/\frak b$, we have $\dim \overline{R}\leq 3$ and further $M$ and $N$ are $R/\frak b$-module. Clearly $(0:_{\Coker f}\frak a\overline{R})$ is a finitely generated $\overline{R}$-module. It follows from [NS, Theorem 2.8] that 
$\ker f$ and $\Coker f$ are $\frak a\overline{R}$-cofinite and so using \cref{del}, they are $\frak a$-cofinite.
\end{proof}



\end{document}